\documentclass{amsart}
\usepackage{amssymb,graphicx,color}
\usepackage[colorlinks=false,pdfborder={0 0 0}]{hyperref}

\newcommand{\Real}{{ \text{Re}}}
\newcommand{\Ima}{{ \text{Im}}}
\newcommand{\R}{{\mathbb R}}
\newcommand{\C}{{\mathbb C}}
\newcommand{\T}{{\mathbb T}}
\newcommand{\RS}{{\mathbb C_{\infty}}}
\newcommand{\D}{{\mathbb D}}
\newcommand{\capacity}{{\operatorname{cap}}}
\newcommand{\dist}{{\operatorname{dist}}}

\newtheorem{thm}{Theorem}
\newtheorem{prop}{Proposition}

\newtheorem{lemma}{Lemma}

\theoremstyle{definition}
\newtheorem{prob}{Problem}
\newtheorem{defini}{Definition}

\theoremstyle{remark}
\newtheorem{remark}{Remark}

\title{The P\'olya-Tchebotar\"ov problem}
\author{Joaquim Ortega-Cerd\`a}
\address{Dept.\ Matem\`atica Aplicada i An\`alisi,
 Universitat  de Barcelona,
Gran Via 585, 08071 Bar\-ce\-lo\-na, Spain}
\email{jortega@ub.edu}
\author{Bharti Pridhnani}
\address{Dept.\ Matem\`atica Aplicada i An\`alisi,
 Universitat  de Barcelona,
Gran Via 585, 08071 Bar\-ce\-lo\-na, Spain}
\email{bharti@maia.ub.es}
\thanks{Supported by DGICYT grant MTM2008-05561-C02-01 and the CIRIT grant
2005SGR00611}
\date{\today}
\subjclass[2000]{Primary }
\begin{document}

\begin{abstract}
We describe the solutions to the problem of identifying the continuum in the
complex plane that minimizes the logarithmic capacity among all the continuum
that contain a prefixed finite set of points. This description can be
implemented numerically and this can be used to improve the estimates on the
Bloch-Landau constant and other related problems as the maximal expected
lifetime of the Brownian motion on domains of inner radius one or the principal
eigenvalue for the Laplace operator on such domains.
\end{abstract}
\maketitle

\section{Introduction and history of the problem}
P\'olya in \cite{Polya29} discussed the following problem which was suggested
to him by Tchebotar\"ov:
\begin{prob}\label{cap}
Given a finite number of points $E:=\left\{a_{1},\ldots,a_{n}\right\}\subset
\C$, find the continuum $K$  with minimal logarithmic
capacity such that $E\subset K$.
\end{prob}
For any continuum $K$, its complement in the Riemann Sphere,
$\Omega=\RS$ is simply connected, therefore there exists a unique conformal map
$f:\Omega\rightarrow\RS\setminus D(0,R)$ such that
$f(z)=z+c_{0}+c_{-1}/z+\cdots$ around $\infty$. Here $R=R(\Omega,\infty)$ is
called the conformal radius of $\Omega$ with respect to $\infty$. In fact,
$\capacity(K)=R(\Omega,\infty)$. This provides an equivalent formulation of
Problem~\ref{cap}, which is usually called the outer reformulation of the
P\'olya-Tchebotar\"ov problem:
\begin{prob}\label{Laur}
Given a finite number of points $E:=\left\{a_{1}, \ldots,a_{n}\right\}\subset
\C\setminus\left\{0\right\}$ find a conformal map $f:\D\to \C\setminus
E$ such that $f(0)=0$ and $\left|f'(0)\right|$ is maximal.
\end{prob}
So we are looking for a simply connected domain $\Omega$ that contains the
origin, it is contained in $\C\setminus E$ and such that the density of the
hyperbolic metric at the origin is minimal. Such domain will be called an
extremal domain and the corresponding conformal map, an extremal map.

The existence of the solution is obvious by a normal family argument. This
problem was studied in detail by Laurentiev in \cite{Laurentiev30}. He proved
the uniqueness and the basic structure of the solution by the method of
variations of the boundary. The structure of the extremal domain is
characterized by the following theorem, see \cite{Laurentiev34}.
\begin{thm}[Laurentiev]\label{Laurentiev}
Given a finite number of points
$E=\left\{a_{1},\ldots,a_{n}\right\}\subset \C$, there exists a unique
extremal domain $\Omega=f(\D)$ for the problem~\ref{Laur} and it is
characterized by the
following properties:
\begin{enumerate}
	\item Each point of the plane belongs to either $\Omega$ or
$\Gamma:=\partial \Omega$.
	\item The boundary $\Gamma$ consists of finitely many simple arcs of
analytic curves. The points $a_{i}$ and $\{\infty\}$ are endpoints of $n+1$
distinct arcs. Every point of $\Gamma$ different from the $a_{i}$ or
$\{\infty\}$ either belongs to a unique arc and it is a regular point of
$\Gamma$, or it is the common end of at least three arcs.
	\item\label{hs} To any arc $\alpha\beta$ consisting of regular points of
$\Gamma$ there correspond under the conformal mapping $f^{-1}$ two arcs of the
same length on the unit circle.
\end{enumerate}
\end{thm}
When the last property \ref{hs} holds we say that the arcs are
\emph{harmonically symmetric} with respect to the origin and it will be the key
property to find a numerical algorithm to determine the solution to the problem
mentioned above.

In the proof of this last theorem, Laurentiev assumed that the desired domain is
bounded by finitely many simple Jordan arcs. This assumption was removed by
Goluzin who used the method of inner variations to prove the following:
\begin{thm}[Goluzin, \cite{Goluzin46}]
Let $a_{1}, \ldots,a_{n}$ be arbitrary given points in $\C$. Let
$K$ be the extremal continuum for
Problem~\ref{cap}.
Then $K$ is the union of the closures of all critical trajectories of the
quadratic differential
\[
Q(z)dz^{2}=-\frac{\prod^{n-2}_{l=1}(z-b_{l})}{\prod^{n}_{k=1}(z-a_{k})}dz^2
\]
where $b_{l}$ are some unknown parameters. 
The extremal univalent function
$g:\RS\setminus \D \to \RS\setminus\{a_1,\ldots,a_n\}$, with $g(\infty)=\infty$
that maximizes $g'(\infty)$ must satisfy the following differential equation
\[
\left(zg'(z)\right)^{2}=\frac{\prod^{n}_{i=1}(g(z)-a_{i})}{
\prod^{n-2}_{j=1}(g(z)-b_{j})}.
\]
\end{thm}
\begin{remark} The points $b_i$ correspond to common end points of several
arcs. If some point $b_{i}$ is a common end of $m$ arcs, then the
term $g(z)-b_{i}$ will appear exactly $m-2$ times in the differential equation.
\end{remark}
Later we will explain how to use this differential equation to obtain a
numerical solution to Problem~\ref{Laur}.

Goluzin gave a more general result where the problem is to maximize
$|f^{(n)}(0)|$ for any $n\ge 1$. An account of his work is in
\cite[Chap.~4]{Goluzin69}. By using this description and after considerable
work, Kuzmina in \cite{Kuzmina82} computes the extremal domain in the case of
three points and in  \cite{Fedorov84} this is extended to four points with a
certain symmetry (two of the points must be symmetric with respect to a line
that
passes by the other two).

Later on, Tamrazov found an explicit solution for the problem of $n$ points. The
general solution to Problem~\ref{Laur} is, according to \cite{Tamrazov05},  of
the form
\[
 f(z)=\int^z (\zeta-1)^{-3}\Bigl(\prod_{\alpha\in V\setminus\{1\}}
(\zeta-\alpha) \Bigr)\prod_{v\in W(\Gamma)}
\Bigl(\prod_{\beta\in W_v}(\zeta -\beta) \Bigr)^{1/\tau(v)}\, d\zeta,
\]
where $\alpha$ and $\beta$ are a finite number of points in $\T$ and $\tau(v)$
is a positive integer. This result does not seem to be completely clear, because
the function $f$ corresponds to a Schwartz-Christoffel formula, thus $f(\T)$ is
going to be a collection of straight segments (one of them going to $\infty$)
but, even when we have only three points, in most cases (except of very
symmetric ones) the solution to the P\'olya-Tchebotar\"ov problem are not
straight lines.

Nevertheless the main idea of Tamrazov paper, that all solutions can be exactly
parametrized by a finite planar graph is indeed correct. We will give a
different proof of this fact. Our approach although it will not yield an
``explicit'' formula it will be constructive and it is possible to implement a
numerical algorithm that produces an approximation to the solution of the
P\'olya-Tchebotar\"ov problem.

On Sections~\ref{sec2} and \ref{sec3} we prove that all solutions are codified
by a ``nested partition'' which are defined there. This is more convenient for
us, although it is completely equivalent to a parametrization by graphs. To each
set of points the extremal continuum is in correspondence with  a unique
``nested partition'' of $\T$ and conversely, each ``nested partition'' provides
a solution.  Thus all the combinatorial data of the solution is codified in
these partitions.

Once we have a parametrization of all possible solutions we introduce in
\ref{sec4} a numerical algorithm to compute numerically the solutions (i.e. to
determine the parameters) to the P\'olya-Tchebotar\"ov problem. We illustrate
the method making it explicit in the case of $3$ points and $6$ points (with a
certain symmetry). This last case is particularly interesting because it will be
of use for the applications that we had in mind which are developed in
Section~\ref{sec5}.

The P\'olya-Tchebotar\"ov problem is rather basic, thus it is not
surprising that it arises in connection with many other problems. The most
evident case is in the estimates of the univalent Bloch-Landau constant, the
precise formulation of the problem is in Section~\ref{sec5}. This has been
exploited in \cite{CaOr08} where this constant was improved. This work is its
natural continuation. Here we will provide more
sophisticated examples  and we will use the same type of
domains to improve the estimates of two other extremal problems that were
introduced in \cite{BaCa96}: the expected lifetime of the Brownian motion in a
domain with inner radius one and the estimation of the principal frequency of
such domains. For the precise definitions and results, see again
Section~\ref{sec5}. 

There are other potential applications of the
P\'olya-Tchebotar\"ov problem which could benefit from our (numerical) solution.
For instance the best estimates in the Smale mean value conjecture obtained in
\cite{Crane07} rely on the computation of the solution of the problem with
three points. We have not pursued improvements on this problem.

\textbf{Acknowledgment: } We are indebted to \`Alex Haro for illuminating
conversations about the numerical implementation of the algorithm.

\section{The parametrization of all solutions}\label{sec2}
In view of Laurentiev and Goluzin results ,the continuum that is a solution 
of  Problem~\ref{cap} form a finite planar tree with endpoints in the points
$\{z_1,\ldots,z_n\}$. The remaining nodes are of order at least three. This
motivates the following definition introduced in \cite{Tamrazov05}:

\begin{defini}
A graph $\Gamma$ is a P\'olya-Tchebotar\"ov graph (PT-graph for short) if it is
a finite planar tree with the properties: 
\begin{enumerate}
 \item All sides of the graph are linear segments.
 \item There are no nodes of order 2.
 \item The sum of the length of all sides is exactly 1/2.
\end{enumerate}
\end{defini}
We say that $\Gamma$ is normalized if we  mark one of the vertex of the tree (a
node of order one). 
Two PT-graphs are equivalent if there is an isometry from one to the
other that extends to an orientable homeomorphism of the plane. If they are
normalized we also require that the isometry sends the marked vertex from the
first graph to the marked vertex of the second graph. We will talk of a
PT-graph $\Gamma$ to denote the whole equivalence class.

Since Problem~\ref{cap} is invariant by translation we can be normalize the
data $\{z_1,\ldots,z_n\}$ to assume that $z_1=0$. 

The main result is that there is a natural way of parametrizing the solutions
to the P\'olya-Tchebotar\"ov problem by normalized PT-graphs. That is for
any graph there
is associated a  unique continuum that solves a P\'olya-Tchebotar\"ov problem
and conversely
all solutions arise in this way. 

\begin{figure}[ht]
\includegraphics[width=6cm]{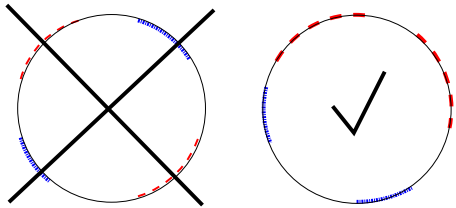}
\caption{Not nested and nested pairs}\label{twin}
\end{figure}

Let us describe how to associate a solution to each graph. First we need one
further definition.
\begin{defini}
A partition $\Pi$ of the unit circle in a finite number of intervals is called a
properly nested partition if the following properties hold:
\begin{enumerate}
 \item The intervals come in pairs of equal length, i.e.
\[
 \Pi=\{I_1,\ldots,I_n\}\cup\{J_1,\ldots,J_n\}
\]
and $|J_i|=|I_i|$ for all $i=1,\ldots,n$.
\item There are no nested pairs of intervals, i.e. a couple $I_i,J_i$
never separates another couple $I_k,J_k$. See Figure~\ref{twin}
\end{enumerate}
\end{defini}
Since all pairs are not nested we can be sure that there exist at least
one pair $I_i,J_i$ of adjacent intervals. Two partitions are equivalent if one
rotation sends one to the other. They are normalized if we mark one of the
adjacent pairs of intervals.

It is easy to see that any PT-graph provides a nested partition and
conversely.
We start from a vertex of the graph and travel through its edges directwise. To
each edge of the tree we associate an interval in the circle of the same length
(on the unit circle we consider the normalized length). Each edge of the tree is
visited twice, once for every side. We consider the pairs of intervals $I_i,J_i$
to intervals that correspond to  different sides of the same edge of the tree.

\begin{figure}[ht]
\includegraphics[width=3cm]{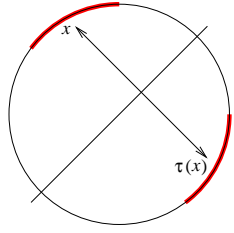}
\caption{The involution $\tau$ restricted to a pair of
intervals}\label{involution}
\end{figure}

Finally for any given properly nested partition there is associated an
involution $\tau$ defined on the circle (except in a finite number of points
corresponding to the end points of the intervals). Two points $x,y$ are related
by the involution $\tau$ if $x$ belongs to the interval $I_i$ and $y$ belongs to
its pair $J_i$. The definition of $\tau$ in each of the intervals  is the
reflection along the diameter of the disk that passes halfway in between the
pair of intervals as in Figure~\ref{involution}

We are going to prove a ``welding'' type theorem:
\begin{thm}
For any given properly nested partition $\Pi$ of $\T$ and its associated
involution $\tau$
there is a conformal map $f_\Pi:\D\to \C\setminus\Gamma$ where $\Gamma$ is a 
finite union of analytic arcs $\Gamma=\cup_{i=1}^n \gamma_i$ such that the 
image of any pair $I_k,J_k\subset \Pi$  is one arc
$\gamma_k$ and $f_\Pi(x)=f_\Pi(\tau(x))$ for all $x$ in the intervals. Moreover
$f_\Pi$ is unique up to postcomposition with automorphisms of $\C$.
\end{thm}
\begin{remark}
If we compose $f_\Pi$ with a translation we obtain a function that satisfies
the Laurentiev conditions of Theorem~\ref{Laurentiev}. Thus for any $\Pi$ we
get a solution to Problem~\ref{Laur}. The converse is even more clear.
Laurentiev theorem shows that the boundary of the extremal domain $\partial
\Omega$ is a tree, that is homeomorphic to a rectilinear PT-graph
$\Gamma$. The length of each edge of $\Gamma$ is one half of the harmonic
measure of the corresponding edge of $\partial\Omega$. The conformal map
$g:\Omega\to D$ with $g(0)=0$ gives a partition $\Pi$ of the unit circle.
By property~\eqref{hs} of Laurentiev theorem, we get $g^{-1}=f_\Pi$.
\end{remark}

\section{Proof of the welding Theorem}\label{sec3}
Given the partition we will proceed to construct the mapping in a finite number
of steps. In each step the following lemma is the key
\begin{lemma}\label{lemabasic}
Given two adjacent intervals $I,J\subset \T$ in the circle, such that $\T\neq
\{I\cup J\}$ and a quasisymmetric homeomorphism $\tau: I\to J$ that
fixes the common point $p$, there is a simple Jordan arc $\gamma\subset \D$ with
one endpoint at $1$ and a conformal mapping $f:\D\to \D\setminus \gamma$ such
that $f(0)=0$ and $f(x)=f(\tau(x))\in\gamma$ for all $x\in I$. The mapping $f$
and the curve $\gamma$ are unique, they depend on $I,J$ and $\tau$.
\end{lemma}

\begin{proof}
Let $I_0=\{e^{it};\ t\in[0,\pi/2]\}$ and $J_0=\{e^{it};\ t\in[-\pi/2,0]\}$,
and let $\tau_0:I_0\to J_0$ be defined as $\tau_0(z)=\bar z$. If these were the
data of the problem, it will be readily solved by the mapping $f_0(z)=...$ that
maps $\D$ to $\D\setminus [r,1)$. In the general situation, there exists an
homeomorphism of the circle $T: \T\to\T$ such that $T(I)=I_0$, $T(J)=J_0$ and 
$\tau_0(T(x))=T(\tau(x))$ for all $x\in I$. In $I$, the map $T$ is defined
linearly. On $J$ we use as definition $T(y)=\tau_0(T(\tau^{-1}(y)))$ and
outside $I$ and $J$ we define it linearly. Since $\tau$ is asymmetric then
$T$ is quasisymmetric. In general $T$ has the same regularity as $\tau$. By the
Beurling-Ahlfors extension theorem \cite{BeurAhl56} there is a quasiconformal
homeomorphism of the
disk that extends $T$. We will still denote it by
$T$. 

Let $g=f_0\circ T$. The map $g$ is mapping $I$ and $J$ to the arc $[r,1]$ in
such a way that $g(x)=g(\tau(x))$ for all $x\in I$. It is not a conformal, but
it is a quasiconformal map because $T$ is quasiconformal and $f_0$ is conformal.
This can be corrected by solving a Laplace-Beltrami equation. We want an
homeomorphism $\Psi:\D\to\D$ such that $k=\Psi\circ g$ is conformal. Thus
$k_{\bar z}= \Psi_{\bar z}\overline{g_z} +\Psi_{z} g_{\bar z}=0$. There is
always a solution $\Psi$ to this Laplace-Beltrami equation that is an
homeomorphism from $\D$ to $\D$ since $|g_{\bar z}/g_z|<k<1$ almost everywhere
by the measurable Riemann mapping theorem. Finally we compose $k$ with an
automorphism of the disk $\kappa$ and the desired function is $f=\kappa\circ
k$. The automorphism $\kappa$ is chosen to make sure that $f(0)=0$ and 
the endpoint of $f(I)=\Psi([r,1])$ is at $1$. The regularity of $f$ at the
boundary is as good as that of $\Psi$ and that itself is determined by the
regularity of $\tau$.

If there were two curves $\gamma_0$ and $\gamma_1$ and two maps $f_0$ and
$f_2$ with the same property then $g=f_1\circ f_0^{-1}$ is a conformal mapping
from  $\D\setminus \gamma_0$ to $\D\setminus\gamma_1$ such that $g(0)=0$.
Moreover $g$ extends continuously to $\gamma_0$ because for any point in
$\gamma_0$ the preimage by $f_0$ are two points $\{x,y\}$ in $\T$ that are
related by $\tau$, thus $f_1(x)=f_1(y)$, therefore $g$ extends continuously to
$\gamma_0$ and therefore it extends analytically, thus $g(z)=e^{i\theta} z$.
Since both $\gamma_0$ and $\gamma_1$ start have an endpoint in $1$, then
$g(z)=z$ and $f_0=f_1$. 
\end{proof}
With the same proof, mutatis mutandi, we deal with the case $\T=I\cup J$ and we
obtain
\begin{lemma}\label{variant}
Given two adjacent intervals $I,J\subset \T$ in the circle, such that $\T=
\{I\cup J\}$ and a quasisymmetric homeomorphism $\tau: I\to J$ that fixes the
common points, there is a simple Jordan arc $\gamma\subset \C$ with one endpoint
at $1$ and the other at $\infty$ and a conformal mapping $f:\D\to \C\setminus
\gamma$ such that $f(0)=0$, $f'(0)>0$ and $f(x)=f(\tau(x))\in\gamma$ for all
$x\in I$. The mapping $f$ and the curve $\gamma$ are unique, they depend on
$I,J$ and $\tau$.
\end{lemma}

\begin{proof}[Proof of the theorem]
We take any pair of adjacent intervals $(I_k,J_k)$ in the
partition $\Pi$ corresponding by the involution $\tau$. There are always
adjacent pairs because the partition is properly nested (they
correspond to edges with and endpoint in a vertex of the graph). Applying the
Lemma we find a conformal mapping $f_1$ that welds together the pair of
intervals in a curve $\gamma_1\subset \D$. The mapping $f_1$ induces a new
partion $\Pi_1$ of $\T$, $I_j^1=f_1(I_j)$ and $J_j^1=f_1(J_j)$ for all $j$
except for the pair $(I_k,J_k)$ that was welded together. This new partition is
again correctly nested because the order is preserved except for a pair of
adjacent intervals that ``collapses''. The number of pairs of intervals is
one less than in $\Pi_1$. The intervals in each pair are no
longer of the same size but nevertheless the map $f_1$ induces a new involution
on them, $\tau_1=f_1\circ\tau\circ f_1^{(-1)}$. 

Now we repeat the procedure. We take any other new pair of adjacent
of intervals of the new partition $\Pi_1$ corresponding by $\tau_1$ and we glue
them together. This can be done by a mapping $f_2$ applying the Lemma because
$\tau_1$ is quasisymmetric (it is in fact piecewise real analytic). In this way
we get again a new involution $\tau_2=f_2\circ\tau_1\circ f_2^{(-1)}$ and a new
nested partition $\Pi_2$. 

In this way we keep gluing pairs of intervals until we are left only with two
intervals and an involution $\tau_n$ that relates them. In this last step we
use Lemma~\ref{variant} to get $f_n$. The final conformal mapping is
$f_\Pi=f_n\circ\cdots \circ f_1$.

There is basically only one such map $f_\Pi$ (except for composition with maps
of the form $az+b$). The proof is as in Lemma~\ref{lemabasic}, suppose there is
another such map $g_\Pi$. Let $\Gamma=f_\Pi(\T)$, then $h=g_\Pi\circ
f_\Pi^{-1}$ is a one to one analytic mapping  $h:\C\setminus\Gamma\to \C$.
Moreover since the preimage of any regular point in $\Gamma$ are two points in
the circle that are related by the involution, and $g$ maps both points to the
same point, then $h$ can be extended continuously to a conformal map from $\C$
to $\C$, thus $h(z)=az+b$. 
\end{proof}

\section{Numerical Algorithm to find solutions}\label{sec4}
As we mentioned above, in order to implement numerically an algorithm to find
the solution to Problem~\ref{Laur}, we used an important property of the
extremal domain and the differential equation obtained by Laurentiev. Denote by
$\Omega_{n}$ the desired extremal domain for the Problem~\ref{Laur} in case of
$n+1$ points $\left\{a_{1}, \ldots,a_{n},\infty\right\}$. Let
$f:\D\rightarrow\Omega_{n}$ be the conformal map such that $f(0)=0$. We know
that $f$ satisfies the following differential equation
\begin{equation}\label{eqgen}
\left(\frac{zf'(z)}{f(z)}\right)^{2}=C\frac{\prod^{n}_{i=1}(f(z)-a_{i})}{\prod^{
n-1}_{j=1}(f(z)-b_{j})},
\end{equation}
where the parameters $b_{j}$ are unknown and
$C=\frac{\prod^{n-1}_{l=1}(-b_{l})}{\prod^{n}_{k=1}(-a_{k})}$. Using the
solutions of this
differential equation and the last property of Theorem~\ref{Laurentiev} we have
implemented the resolution of Problem~\ref{Laur} for some cases of $n$. The
system becomes more delicate as $n$ increases, the combinatorics and the
dimensions of the systems to solve become bigger. We will show in detail the
solution in the case $n=3$ to illustrate the method and $n=6$ with some extra
symmetry, because this will be enough for the applications that we have
in mind. The code where this algorithm is implemented (for four points and 6
points with symmetry) can be downloaded from
\href{http://www.maia.ub.es/cag/code/tchebotarev/}
{http://www.maia.ub.es/cag/code/tchebotarev/}.

\subsection{Case of 3 points}
\begin{figure}[ht]
    \begin{center}
    \includegraphics[width=7cm]{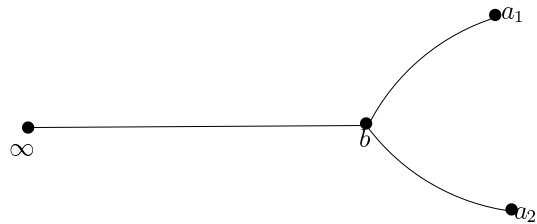}
    \end{center}
    \caption{Sketch of the extremal compact for three points}{\label{3points}}
\end{figure}
\begin{figure}[ht]
    \begin{center}
    \includegraphics[width=10cm]{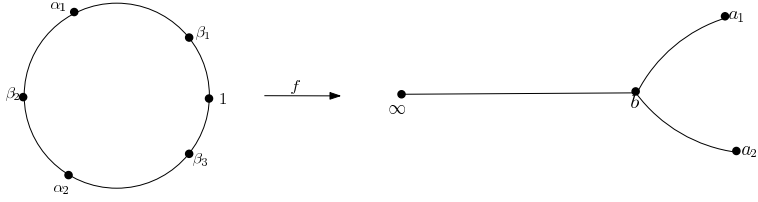}
    \end{center}
    \caption{Configuration for $n=3$}{\label{3pointsconfig}}
\end{figure}

Let's start with 3 points. Assume that we have $a_{1},a_{2},a_{3}$ three points
such that $a_{3}=\infty$ and $a_{1},a_{2}\neq 0$. Without loss of generality we
will always assume that $f(1)=\infty$. In the case of three points the extremal
domain is very clear
(see Figure~\ref{3points}). We only have one unknown parameter denote it by $b$
in the differential equation \eqref{eqgen} that reduces to:
\begin{equation}\label{eq3p}
f'(z)^{2}=C\frac{(f(z)-a_{1})(f(z)-a_{2})}{f(z)-b}\frac{f(z)^{2}}{z^{2}}
\end{equation}
Recall that to any regular arc $\alpha\beta$ of $\partial\Omega$ there
corresponds two arcs with equal lengths on the unit circle. So this means that
we have the configuration on the unit circle shown in
Figure~\ref{3pointsconfig},
where $f(0)=0$, the arcs $(1\beta_{1}),(\beta_{3}1)$ are mapped to the arc
$\infty b$, $\beta_{1}\alpha_{1}$ and $\alpha_{1}\beta_{2}$ are mapped to the
arc $ba_{1}$ and the arcs $\beta_{2}\alpha_{2},\alpha_{2}\beta_{3}$ into
$ba_{2}$. Note that $f(e^{i\alpha_{i}})=a_{i}$ for $i=1,2$ and
$f(e^{i\beta_{i}})=b$ for $i=1,2,3$.

The solution of the problem can be viewed as the solution of a system of
non-linear equations. If we know the value of $f'(0)$ and $b$ we can compute the
coefficients of the mapping $f$ using the differential equation
\eqref{eq3p}. We computed also the values of $\alpha_{i},\beta_{i}$. Note
that as the arcs $1\beta_{1}$,$\beta_{3}1$ must have same length, we get
$\beta_{3}=2\pi-\beta_{1}$ (we will always take the angles in the range
$[0,2\pi)$). So, we have 6 real unknown parameters in our problem:
$\Real(f'(0)),\Ima(f'(0)),\Real(b),\Ima(b),\beta_{1},\beta_{2}$ and using the
last property of Theorem~\ref{Laurentiev} we can impose the following three
complex equations
\[
\begin{cases}
	f(e^{i\beta_{1}/2})=f(e^{-i\beta_{1}/2}).\\
	f(e^{i(\alpha_{1}+\beta_{1})/2})=f(e^{i(\alpha_{1}+\beta_{2})/2}).\\
	f(e^{i(\alpha_{2}+\beta_{2})/2})=f(e^{i(\alpha_{2}+\beta_{3})/2}).
\end{cases}
\]
We used a the hybrid method to find an approximation of the roots of the system
(see \cite{Powell70} for more details of the method).

To apply the root-finding method we need to evaluate $f(e^{i\gamma})$ for any
$\gamma\in[0,2\pi)\setminus\left\{\alpha_{1},\alpha_{2},\beta_{1},\beta_{2},
\beta_{3}\right\}$. For that, denote $z(t)=f(te^{i\gamma})$. We know that
$z(0)=0$ and $z'(0)=f'(0)e^{i\alpha}$. Note that $z(1)=f(e^{i\alpha})$. We can
get the differential equation satisfied by $z(t)$ and solve it to obtain the
value in time $t=1$. We get 
\begin{equation}\label{edo}
z'(t)^2=C\frac{(z(t)-a_{1})(z(t)-a_{2})}{(z(t)-b)}\frac{z(t)^2}{t^2}
\end{equation}
Note that this equation only defines $z'(t)$ up to a sign, we will deal with
this problem by analytic continuation. Once we fix the derivative at the origin
there is a single analytic branch that solves the equation. To solve it we
used the Taylor integration method which allows us to integrate the singularity
in $t=0$. As $f$ is conformal, we know that $z(t)=z_{1}t+z_{2}t^{2}+\ldots$,
where
$z_{1}=f'(0)e^{i\gamma}$. Now if we do the calculations in the equation
\eqref{edo} we get a recurrence for the coefficients till the order we want.
Now we can estimate the radius of convergence of the
obtained series. And therefore proceed using Taylor method to integrate the 
differential equation until $t=1$. Hence we will be able to impose the equations
to solve our problem.
\subsection{Case of 6 points with symmetry}\label{sub6p}
\begin{figure}[ht]
    \begin{center}
    \includegraphics[width=12cm]{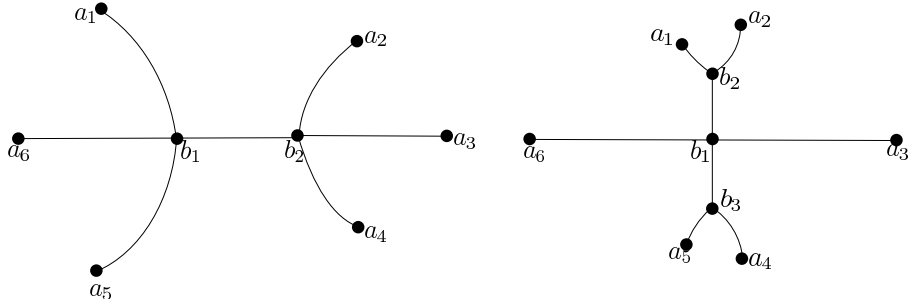}
    \end{center}
    \caption{Structure of the extremal domains for $n=6$ with
symmetry}{\label{6pointsymm}}
\end{figure}

Now consider the case of $6$ points $a_{1},a_{2}, \ldots,a_{6}$ such that
$a_{6}=\infty$, $a_{3}\in\R$ and $a_{5}=\bar{a_{1}},a_{4}=\bar{a_{2}}$. The
extremal compact in this case may be of two types (see
Figures~\ref{6pointsymm}, \ref{6points1config} and \ref{6points2config}).
So, we have two type of configurations in this special case
\begin{enumerate}
	\item $0\ \beta^{1}_{1}\ \alpha_{1}\ \beta^{1}_{2}\ \beta^{2}_{1}\
\alpha_{2}\ \beta^{2}_{2}\ \alpha_{3}\ \beta^{2}_{3}\ \alpha_{4}\ \beta^{2}_{4}\
\beta^{1}_{3}\ \alpha_{5}\ \beta^{1}_{4}\ 2\pi$
	\item $0\ \beta^{1}_{1}\ \beta^{2}_{1}\ \alpha_{1}\ \beta^{2}_{2}\
\alpha_{2}\ \beta^{2}_{3}\ \beta^{1}_{2}\ \alpha_{3}\ \beta^{1}_{3}\
\beta^{3}_{1}\ \alpha_{4}\ \beta^{3}_{2}\ \alpha_{5}\ \beta^{3}_{3}\
\beta^{1}_{4}\ 2\pi$,
\end{enumerate}
where $f(e^{i\beta^{j}_{k}})=b_{j}$ for $j=1,2$ and $k=1,2,3,4$,
$f(e^{i\alpha_{i}})=a_{i}$ for $i=1,2,3,4,5$ and $f(1)=\infty$.
Using the symmetry we can do some reductions to get a system of equation with
less dimension. For example, we
can always assume that the point $a_{3}\in\R$, so that $b_{1},b_{2}\in\R$ for
the first configuration and $b_{1}\in\R$ for the second one. Note that in this
last configuration, by the symmetry of the problem, $b_{3}=\bar{b_{2}}$.
Moreover as this is a symmetric case, $f'(0)$ must be real and
$\alpha_{3}=\pi$.

\subsubsection*{First configuration}
In this case we have $7$ real unknown parameters: $\Real(f'(0))$,
$\Real(b_{1})$,
$\Real(b_{2})$, $\beta^{1}_{1}$, $\beta^{1}_{2}$, $\beta^{2}_{1}$,
$\beta^{2}_{1}$ and we can impose the following equations
\[
\begin{cases}
	\Ima(f(e^{i\beta^{1}_{1}/2.0}))=0\\
	f(e^{i(\alpha_{1}+\beta^{1}_{1})/2})=
              f(e^{i(\alpha_{1}+\beta^{1}_{2})/2})\\
	\Ima(f(e^{i(\beta^{1}_{2}+\beta^{2}_{1})/2}))=0\\
	f(e^{i(\alpha_{2}+\beta^{2}_{1})/2})=
              f(e^{i(\alpha_{2}+\beta^{2}_{2})/2})\\
	\Ima(f(e^{i(\alpha_{3}+\beta^{2}_{2})/2}))=0
\end{cases}
\]
\subsubsection*{Second configuration}
\begin{figure}[ht]
    \begin{center}
    \includegraphics[width=12cm]{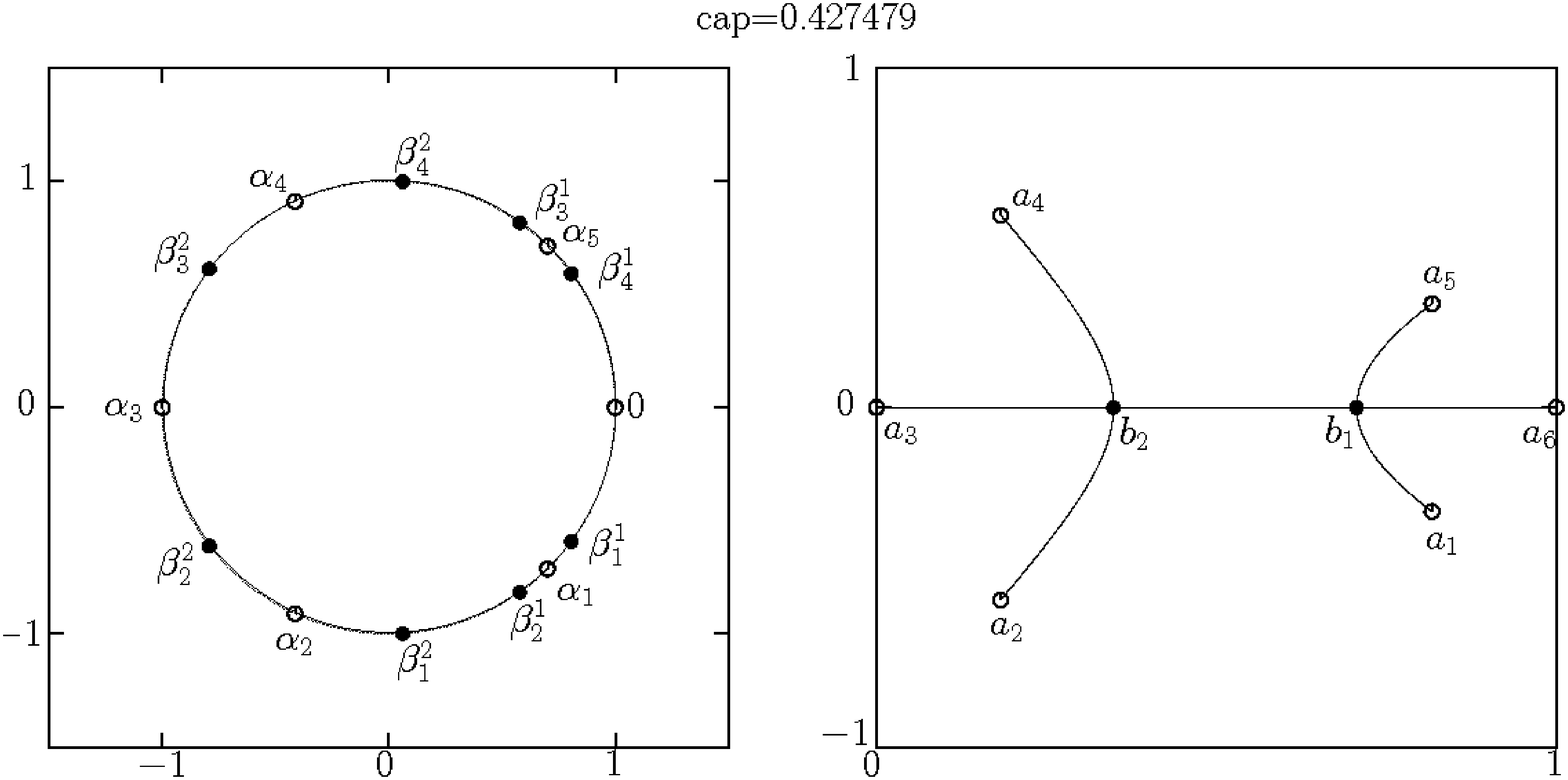}
    \end{center}
    \caption{Extremal domain for $n=6$ with symmetry (configuration
1)}{\label{6points1config}}
\end{figure}

\begin{figure}[ht]
    \begin{center}
    \includegraphics[width=12cm]{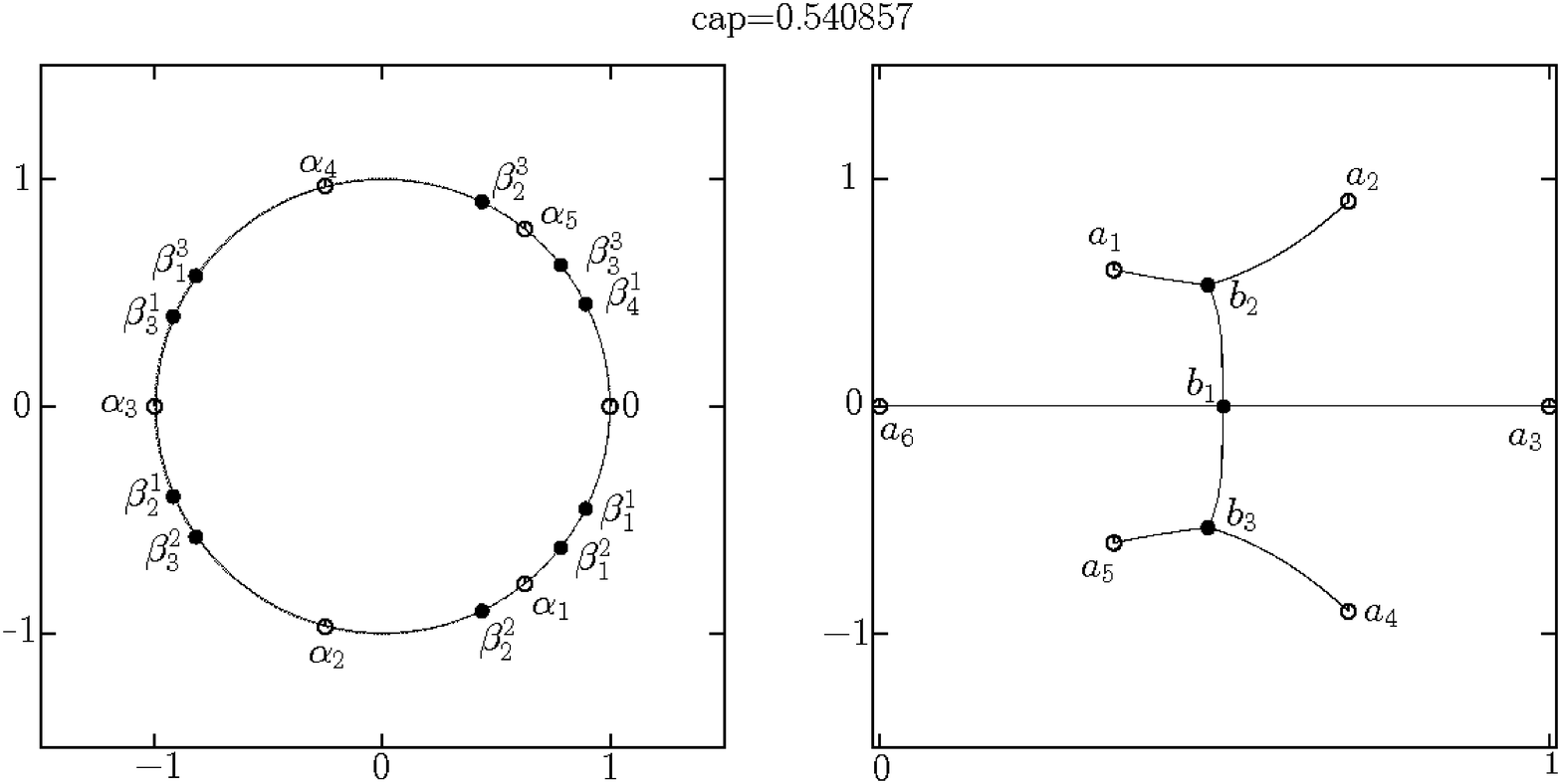}
    \end{center}
    \caption{Extremal domain for $n=6$ with symmetry (configuration
2)}{\label{6points2config}}
\end{figure}

We have $8$ real unknown parameters:
$\Real(f'(0))$, $\Real(b_{1})$, $\Real(b_{2})$, $\Ima(b_{2})$,
$\beta^{1}_{1}$, $\beta^{2}_{1}$, $\beta^{2}_{2}$,$\beta^{2}_{3}$ and we can
impose the
following equations
\[
\begin{cases}
	\Ima(f(e^{i\beta^{1}_{1}/2.0}))=0\\
	f(e^{i(\alpha_{1}+\beta^{2}_{1})/2})=
             f(e^{i(\alpha_{1}+\beta^{2}_{2})/2})\\
	f(e^{i(\alpha_{2}+\beta^{2}_{2})/2})=
            f(e^{i(\alpha_{2}+\beta^{2}_{3})/2})\\
	\Ima(f(e^{i(\alpha_{3}+\beta^{2}_{3})/2}))=0\\
	f(e^{i(\beta^{1}_{1}+\beta^{2}_{1})/2})=
                  f(e^{i(\beta^{2}_{3}+\beta^{1}_{2})/2})\\
\end{cases}
\]
In Figures~\ref{6points1config} and \ref{6points2config} we show an extremal
domain for some 6 points for each configuration. As in the last case, these
figures represent a conformal map $g$ from the complement of $\D$ onto $\Omega$
such that $g(\infty)=\infty$.
\begin{remark}
Note that the solution of the problem depend continuously on the parameters
$a_{i}$, so if we have one solution for some given set points we can do
continuation to reach to any other set of  points (with the same topological
configuration). This has been used and we did the classic continuation (i.e. for
the new set of points we take as a initial condition the solution of the last
set of points).
\end{remark}

\begin{remark}
In the implementation of the method, we found a problem when the distance
between the arcs on the unit circle is very small, we can't integrate properly
the differential equation because we are near the poles $b_{i}$. However this
can be overcome by a change of variables. This has been implemented for the
special case of the first configuration mentioned in the case of 6 points. In
fact, for the application of the P\'olya-Chebotarev we only needed the values of
$f'(0),b_{1},b_{2}$. This data is enough to obtain the expansion in series of
the mapping $f$. This can be done in the following way: given the points
$a_{i}$, we have an initial guess for the unknown values. So, we compute the
critical orbit starting at point $a_{1}$ till the point with imaginary part
equal to $\Ima(a_{1}/2.0)$ and same for the critical orbit starting at $b_{1}$.
The real part of the two points obtained should match if $b_{1}$,$b_{2}$ are the
desired solution. So this is one real equation. Same can be done for the couple
of points $a_{2}$ and $b_{2}$. So we have two real equations. The last equation
can be $\Ima(f(e^{i0.1}))=0$ (this is valid if the points $a_{1},a_{2}$ are a
bit far from the point $a_{6}$).
\end{remark}
\section{Applications of the P\'olya-Tchebotar\"ov problem}\label{sec5}
\subsubsection*{The fundamental frequency of a domain} In 1965, Endre Makai (see
\cite{Makai65}) proved the following theorem solving a
 problem in the study of vibrating membranes raised by P\'olya and Szeg\"o in
their  book \cite{PoSz51}. In 1978, Hayman (\cite{Hayman78}) unaware of it,
reproved the same result.
\begin{thm}
Let $D$ be a simply connected domain in the complex plane. Let $R_{D}$ be the
inradius of $D$, that is, the radius of the largest disc contained in $D$ and
let $\lambda_{D}$ be the first Dirichlet eigenvalue for the Laplacian in $D$.
There is a universal constant $a$ such that
\begin{equation}\label{consta}
\lambda_{D}\geq \frac{a}{R_{D}^{2}}.
\end{equation}
\end{thm}
There have been many efforts to find the best constant $a$ and to identify the
extremal domain for $a$. Makai's proof also shows that the best $a$ satisfies
$1/4\leq a < \pi^{2}/4$. 

The following lemma is useful for giving upper bounds for this constant (see
\cite[Lemma 1.2]{BaCa96} for a proof):
\begin{lemma}
Let $J_{0}(r)$ be the first Bessel function and $j_{0}$ the smallest positive
zero of $J_{0}$. Assume that $D$ is a simply connected domain. Then
\[
\lambda_{D}\leq
j_{0}^{2}\inf_{F}\left\{\frac{1}{\sum^{\infty}_{n=1}|a_{n}|^{2}\delta_{n}}
\right\}
\]
where
\[
\delta_{n}=n^{2}\frac{\int^{1}_{0}J^{2}_{0}(j_{0}r)r^{2n-1}dr}{\int^{1}_{0}J^{
2}_{0}(j_{0}r)rdr}
\]
and the infimum is taken over all conformal mappings
$F(z)=\sum^{\infty}_{n=0}a_{n}z^{n}$ from the unit disc onto $D$.
\end{lemma}
In \cite{BaCa96} Ba\~nuelos and Carroll proved that $0.619<a<2.13$ and provided
examples of domains which are close to the extremal domain. They did this
relating this problem to two other extremal problems:

\subsubsection*{The expected lifetime of a Brownian motion} Let $B_{t}$ be the
Brownian motion in $D$. Let
$\tau_{D}=\inf\left\{t>0:B_{t}\notin D\right\}$ be the first exit time of
$B_{t}$ from $D$. Let us denote by $E_{z}(\tau_{D})$ the expectation of
$\tau_{D}$ under the measure of the Brownian starting at the point $z$ in $D$.
It is known that there is a universal constant $b$ such that, whenever $D$ is a
planar simply connected domain,
\begin{equation}\label{constb}
\sup_{z\in D}E_{z}(\tau_{D})\leq bR_{D}^{2}.
\end{equation}
As before, we want to know the best value of $b$ and the extremal domain for
this last inequality. It is a fact that if $D$ denotes the unit disc of radius
$R_{D}$ then $\sup_{z\in D}E_{z}(\tau_{D})=R_{D}^{2}/2$. It is known that
$1.584<b<3.228$ (see \cite{BaCa96}). In order to give an improved lower bound
for $b$ it is useful to know the following result (see \cite[Lemma 1.1.]{BaCa96}
for
further details):
\begin{lemma}
Suppose that $F(z)=\sum^{\infty}_{n=0}a_{n}z^{n}$ is a conformal mapping from
the unit disc onto a simply connected domain $D$ with $F(0)=z_{0}$. Then
\[
E_{z_{0}}(\tau_{D})=\frac{1}{2}\sum^{\infty}_{n=1}|a_{n}|^{2}.
\]
\end{lemma}

\subsubsection*{The univalent Bloch-Landau constant} 
\begin{figure}[ht]
    \begin{center}
    \includegraphics[width=7cm]{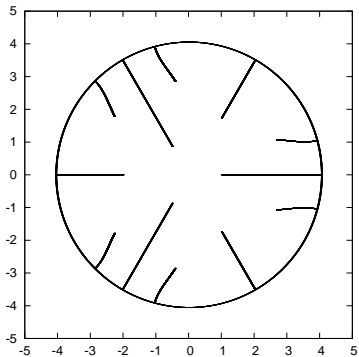}
    \end{center}
    \caption{Structure of the domain in \cite{CaOr08}}\label{domainquim}
\end{figure}

If $f$ is an analytic and
one to one mapping from the unit disc, then there exists a universal constant
$\mathcal U$ such that \begin{equation}\label{constBL}
R_{f(\D)}\geq \mathcal U|f'(0)|.
\end{equation}
This means that the image of the unit disc under any conformal map $f$ contains
discs of radius less that $\mathcal U|f'(0)|$. Note that from Koebe's
1/4-theorem, we know that $\mathcal U\geq 1/4$. The best value of $\mathcal U$
is known as the univalent or schlicht Bloch-Landau constant. We can reformulate
this problem in terms of the density of the hyperbolic metric. If $f$ is a
conformal mapping from the unit disc such that $f(0)=z$ then the density of the
hyperbolic metric is $\sigma(z;D)=1/|f'(0)|$. So we have the following
inequality
\begin{equation}\label{constc}
\sigma_{D}:=\inf_{z\in D}\sigma(z;D)\geq \frac{c}{R_{D}}.
\end{equation}
where $c:=\mathcal U$. From many years there have been efforts to find bounds
for $\mathcal U=c$. This constant was introduced in 1929 by Landau
\cite{Landau29},
who proved that $\mathcal U >0.566$. Reich improved this bound in \cite{Reich56}
($\mathcal U > 0.569$) and Jenkins in \cite{Jenkins61} gave $\mathcal U >
0.57088$.
Many other gave some improved bounds. There are many domains proposed as the
candidate for the extremal domain in order to obtain upper bounds for the
Bloch-Landau constant. For example, Robinson in \cite{Robinson35} proved that
$\mathcal U < 0.658$, Goodman in \cite{Goodman45} that $\mathcal U < 0.65647$
and in \cite{BellHum85} Beller and Hummel proved that $\mathcal U< 0.6564155$.
Finally in \cite{CaOr08} this bound has been improved to $\mathcal U <
0.6563937$. In this last result the resulting domain had all the inner
boundary harmonic symmetric with respect to the origin, see
Figure~\ref{domainquim}.
We will use a slight modification of this last domain, which is still harmonic
symmetric, to give an improved bound of the Bloch-Landau constant. 

As a normalization, we will take domains $D$ with
inradius 1. We will give improved bounds for the constants appearing in the
three problems explained above.

Ba\~nuelos and Carroll in \cite{BaCa96} conjectured that the extremal domain is
the same for all
the three problems. When we restrict these problems to the class of convex
domains this is true. In our work, we will see that a similar domain
improves the bounds for all the three problems.
\subsection{Bloch-Landau constant}
If $h$ is a mapping of the complement of a compact set $E$ onto the complement
of the closed unit disc, it can be expanded (up to a rotation) as
\[
h(z)=\frac{z}{\capacity(E)}+O(1), z\rightarrow \infty.
\]

\begin{figure}[ht]
    \begin{center}
    \includegraphics[width=6cm]{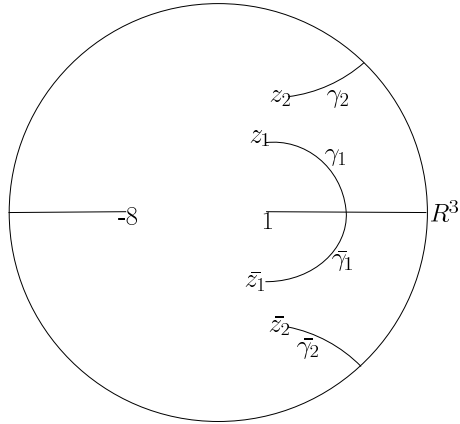}
    \end{center}
    \caption{Domain $\Omega_{z_{1},z_{2},R}$}{\label{domainOmega}}
\end{figure}

\begin{figure}[ht]
    \begin{center}
    \includegraphics[width=6cm]{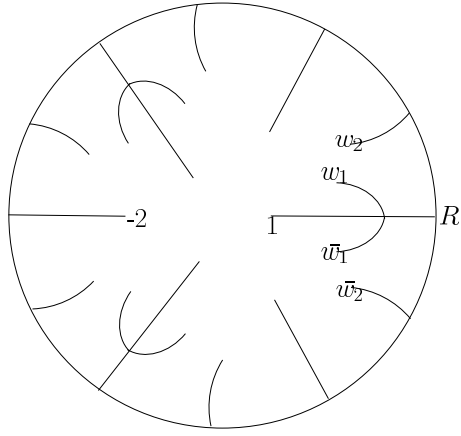}
    \end{center}
    \caption{Domain $D_{\omega_{1},\omega_{2},R}$}{\label{domainDwR}}
\end{figure}

\begin{figure}[ht]
    \begin{center}
    \includegraphics[width=7cm]{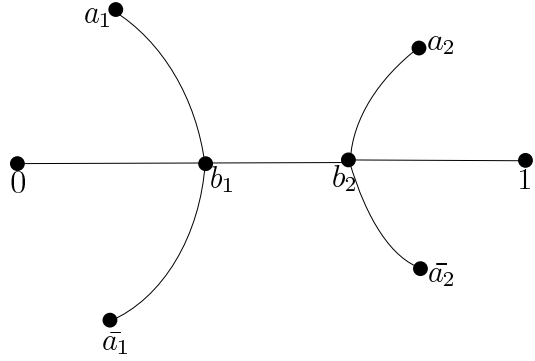}
    \end{center}
    \caption{The continuum E}{\label{E}}
\end{figure}

So we can relate the problem of the extremal domain for the Bloch-Landau
constant with Problem~\ref{Laur} in the case of 6 points because minimizing the
capacity
is equivalent to increasing the derivative at the origin. It is known (see
\cite{Carroll08}) that the arcs making up the extremal configuration must be
harmonically symmetric at infinity. We will work with domains
$\Omega=\Omega_{z_{1},z_{2},R}$ as in Figure~\ref{domainOmega} where $R$ is
bigger than 4 and we chose the arcs
$\gamma_{1},\gamma_{2},\bar{\gamma_{1}},\bar{\gamma_{2}}$ so that this domain is
harmonically symmetric with respect to 0.
If $g$ is a conformal map of $\D$ onto $\Omega_{z_{1},z_{2},R}$ with $g(0)=0$
then $f(z)=z\sqrt[3]{g(z^3)/z^{3}}$ is a conformal map of $\D$ onto the domain
$D_{w_{1},w_{2},R}$ shown in Figure~\ref{domainDwR}. The arcs in this last
domain are harmonically symmetric. We can compute the derivative of $f$ at 0:
$|f'(0)|=\sqrt[3]{|g'(0)|}$. We need this domain to have inradius 1. Later on we
will explain the construction of the domain and the way to get inradius 1.
Let $k(z)=z/(1-z)^{2}$ be the Koebe mapping from the unit disc onto the complex
plane slit along the negative real axis from minus infinity to -1/4.
\begin{prop}
Let $f$ be the conformal map of $\D$ onto $D_{w_{1},w_{2},R}$ such that
$f(0)=0$. Then taking $z_{1}=w_{1}^{3}$ and $z_{2}=w_{2}^{3}$,
\[
\frac{1}{|f'(0)|}=\frac{1}{R}\sqrt[3]{|\psi(-8)-\psi(1)|\capacity(E)}
\]
where $\psi(z)=-1/k(z/R^{3})$ and $E$ is the continuum with minimal capacity
containing 6 given points with symmetry (see Figure~\ref{E}).
\end{prop}
\begin{proof}
Using the above notations, consider the map
\[
\phi(z):=\frac{\psi(z)-\psi(1)}{|\psi(-8)-\psi(1)|}
\]
which maps $\Omega_{z_{1},z_{2},R}$ onto the complement of the continuum $E$
with $a_{1}=\phi(z_{1})$ and $a_{2}=\phi(z_{2})$
($\phi(1)=0,\phi(-8)=1,\phi(0)=\infty$). Note that the harmonic symmetry of the
arcs are preserved since each mapping can be extended continuously to all
internal boundary arcs and the domains involved are symmetric with respect to
the real axis. Let $h$ be the mapping of the complement of $E$ onto the
complement of the unit disc. Then we can define the map $G$ of
$\Omega_{z_{1},z_{2},R}$ onto $\D$ as
\[
G(z):=\frac{1}{h(\phi(z))}.
\]
We can calculate the derivative $G'(0)$ by computing the power series of $G$:
\[
|G'(0)|=\frac{|\psi(-8)-\psi(1)|}{R^{3}h'(\infty)}=\frac{
|\psi(-8)-\psi(1)|\capacity(E)}{R^{3}}.
\]
The capacity of $E$ can be computed numerically as we explained in
\ref{sub6p}. Now
$F:=f^{-1}=z\sqrt[3]{G(z^{3})/z^{3}}$ and the proposition is proved.
\end{proof}

\subsection{Construction of the domain $D_{w_{1},w_{2},R}$}
Now we will explain how to obtain the desired domain in order to have inradius 1
and the derivative $|f'(0)|$ as big as possible. In what follows, $D_{p}$
denotes a disc of radius 1 centered at $p$ and $C_{p}=\partial D_{p}$.

\begin{figure}[ht]
    \begin{center}
    \includegraphics[width=7cm]{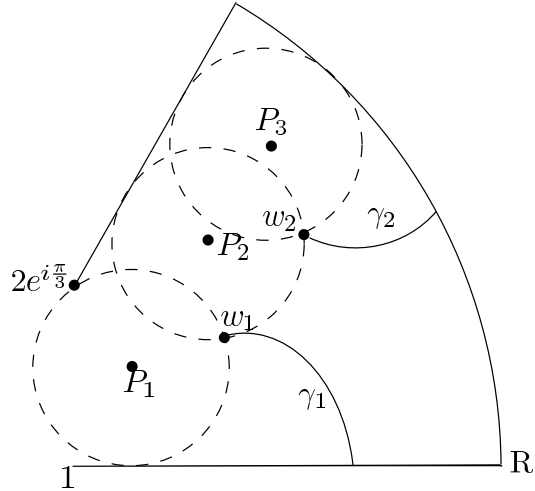}
    \end{center}
    \caption{Election of $w_{1}$ and $w_{2}$}{\label{aplicationdomain1}}
\end{figure}

\begin{figure}[ht]
    \begin{center}
    \includegraphics[width=7cm]{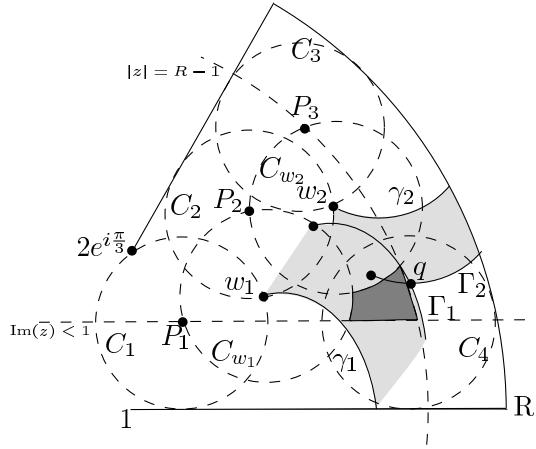}
    \end{center}
    \caption{Prohibited zones}{\label{prohibitedzones}}
\end{figure}

The domain is constructed in stages. Let $D(0,R)$ be a disc centered at the
origin with radius $R>1$. First we remove from this disc three radial slits that
start from the cube roots of the unity. Then, we remove three further slits
starting at two times the cube roots of -1 (these are the first two stages of
Goodman's domain). Now, as we need inradius 1 we need to put some point in order
to not to have discs of radius bigger than 1. Let
$P_{1}=(1+\sqrt{2\sqrt{3}-3},1)$ and denote by $C_{1}$ the circle centered at
$P_{1}$ with radius 1. We have to put some point in $C_{1}$ so that this circle
can't increase. Let $w_{1}$ be a point in this circle. Denote by $C_{2}$ the
circle of radius 1 tangent to the halfline of argument $\pi/3$ containing the
point $w_{1}$. Let $C_{3}$ be the circle of radius 1, tangent to $|z|=R$ and to
the halfline of argument $\pi/3$ and denote by $w_{2}$ the intersection point of
$C_{1}$ and $C_{2}$, $P_{2}$ and $P_{3}$ the centers of $C_{2},C_{3}$,
respectively (see Figure~\ref{aplicationdomain1}).
Now let $\Gamma_{1}$ and $\Gamma_{2}$ be the curves at distance one of
$\gamma_{1}$ and $\gamma_{2}$ (i.e. given a point $\gamma_{i}(t)\in\gamma_{i}$,
let $v_{n}$ be the normalized orthogonal vector to $\gamma_{i}'(t)$, then the
corresponding point in the curve $\Gamma_{i}$ is $\Gamma_{i}(t)=\gamma_{i}(t)\pm
v_{n}$. Let $q$ denote the intersection of $\Gamma_{1}$ and $\Gamma_{2}$. One
sufficient condition to have inradius one is $|q|\geq R-1$. The idea to prove
this is to cover all the points $(x,y)$ so that they can't be centers of circles
(contained in $D_{w_{1},w_{2},R}$) with radius bigger than one. Let's show it
when our points are located in the sector between the segment $[1,R]$ and the
halfline of argument $\pi/3$. If $z=(x,y)$ is a point such that $y\leq
1$,$|z|\geq R-1$ or $\dist(z,T)\leq 1$ (where $T$ is the halfline of argument
$\pi/3$) then obviously we cant have a circle centered at such point with radius
bigger than one. Let $D_{w_{1}}$ and $D_{w_{2}}$ be the discs of radius one
centered at the points $w_{1}$ and $w_{2}$, respectively. If $z\in D_{w_{1}}\cup
D_{w_{2}}$ then $D_{z}$ contains one of the points $w_{1}$ or $w_{2}$.
Therefore, the region $D_{w_{1}}\cup D_{w_{2}}$ is covered. The only risky
region is the one between $\partial D_{w_{1}}$, $\partial D_{w_{2}}$, $|z|=R-1$
and $\left\{z=(x,y)| y=1\right\}$. But this space is covered by the region
delimited by the curves $\Gamma_{1}\cup\Gamma_{2}$ due to the hypothesis (see
Figure~\ref{prohibitedzones}).
\subsection{Results}
We have computed the bounds for all three problems explained before. To
construct the point $w_{1}$ we move on the real axis $x>=1+\sqrt{2\sqrt{3}-3}$
and define the point $P_{2}$ and then $w_{1}$. After this all gets determined.
So given $x$, first we find the biggest $R$ such that $|q|\ge 1-R$,
(because the derivative at the origin will increase with the radius $R$) and
after that we compute the bounds of the constants explained in the three
problems. The results obtained are:

\begin{enumerate}
\item For the Bloch-Landau constant, the best upper bound has been found for
$x=2.1383799965243$ and $R=5.1195152501$ and the improved bound is
\[
\mathcal U\leq 0.656319277272.
\]
The domain obtained is shown in
Figure~$\ref{domainBL}$.
\begin{figure}[ht]
    \begin{center}
    \includegraphics[width=10cm]{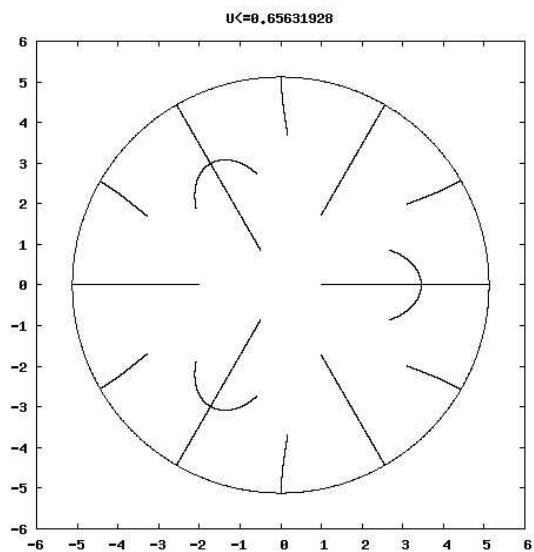}
    \end{center}
    \caption{Domain for the improved upper bound of $\mathcal
U$}{\label{domainBL}}
\end{figure}

\item Computing the coefficients of the conformal mapping obtained for the
domains
$D_{w_{1},w_{2},R}$, the improved upper bound for the fundamental frequency
has been found
for $x=2.1282995811037759$ and $R=5.10223601895443$ and it is
\[
a\leq 2.0907934752309.
\]

\item The improved lower bound for the expected life time of a Brownian
motion has been found for
$x=2.174447128952$ and $R=5.1836816989$ and it is
\[
b\geq 1.670724582110.
\]
\end{enumerate}
In all the computations the error estimates that we required are of the order
of $10^{-12}$, and we are pretty confident on the correctness of the 10 first
digits on the bound, but we have not done a rigorous error analysis.


\def\cprime{$'$}
\providecommand{\bysame}{\leavevmode\hbox to3em{\hrulefill}\thinspace}
\providecommand{\href}[2]{#2}

\end{document}